\documentclass[12pt]{amsart}
\usepackage{amssymb, colordvi}
\usepackage{multirow}
\usepackage{graphicx}

\numberwithin{equation}{section}
\newtheorem{theorem}{Theorem}
\numberwithin{theorem}{section}
\newtheorem{proposition}[theorem]{Proposition}

\newtheorem{corollary}[theorem]{Corollary}

\newtheorem{ex}[theorem]{Example}
\newtheorem{rem}[theorem]{Remark}
\newtheorem{definition}[theorem]{Definition}

\newenvironment{example}{\begin{ex}\rm}{\end{ex}}
\newenvironment{remark}{\begin{rem}\rm}{\end{rem}}
\newcounter{FNC}[page]
\def\fauxfootnote#1{{\addtocounter{FNC}{2}$^\fnsymbol{FNC}$%
     \let\thefootnote\relax\footnotetext{$^\fnsymbol{FNC}$#1}}}

\copyrightinfo{}{}

\renewcommand{\P}{\mathbb{P}}
\newcommand{\C}{\mathbb{C}}
\newcommand{\R}{\mathbb{R}}
\newcommand{\Q}{\mathbb{Q}}
\newcommand{\T}{\mathbb{T}}
\newcommand{\Z}{\mathbb{Z}}

\newcommand{\calA}{\mathcal{A}}
\newcommand{\calB}{\mathcal{B}}
\newcommand{\calW}{\mathcal{W}}

\title{Gale duality for complete intersections} 

\author{Fr\'ed\'eric Bihan}
\address{Laboratoire de Math\'ematiques\\
         Universit\'e de Savoie\\
         73376 Le Bourget-du-Lac Cedex\\
         France}

\email{Frederic.Bihan@univ-savoie.fr}
\urladdr{http://www.lama.univ-savoie.fr/\~{}bihan/}

\author{Frank Sottile}
\address{Department of Mathematics\\
         Texas A\&M University\\
         College Station\\
         Texas \ 77843\\
         USA}
\email{sottile@math.tamu.edu}
\urladdr{http://www.math.tamu.edu/\~{}sottile/}

\thanks{Sottile supported by the Institute for Mathematics and its Applications, 
  NSF CAREER grant DMS-0538734, and Peter Gritzmann of the Technische
         Universit\"at M\"unchen.}  



\begin{document}

\begin{abstract}
  We show that every complete intersection defined by Laurent polynomials in an
  algebraic torus is     
  isomorphic to a complete intersection defined by master functions in the
  complement of a hyperplane arrangement, and vice versa.
  We call systems defining such isomorphic schemes Gale dual systems because the exponents
  of the monomials in the polynomials annihilate the weights of the master functions.
  We use Gale duality to give a Kouchnirenko theorem for the number of solutions
  to a system of master functions and to compute some topological invariants of 
  master function complete intersections.
\end{abstract}
\maketitle

%
\section*{Introduction}
%

A complete intersection with support $\calW$ is a subscheme of the torus
$(\C^\times)^{m+n}$ having pure dimension $m$ that may be defined by a system 
 \[
   f_1(x_1,\dotsc,x_{m+n})\ =\ f_2(x_1,\dotsc,x_{m+n})\ =\ 
    \dotsb\ =\  f_n(x_1,\dotsc,x_{m+n})\ =\ 0
 \]
of Laurent polynomials with support $\calW$.

Let $p_1(y),\dotsc,p_{l+m+n}(y)$ be degree 1 polynomials defining an 
arrangement $\calA$ of hyperplanes in $\C^{l+m}$ and let 
$\beta=(b_1,\dotsc,b_{l+m+n})\in\Z^{l+m+n}$ be a vector of integers. 
A master function of weight $\beta$ is the rational function
\[
   p(y)^\beta\ :=\ 
    p_1(y)^{b_1}\cdot p_2(y)^{b_2}\dotsb p_{l+m+n}(y)^{b_{l+m+n}}\,,
\]
which is defined on the complement $M_\calA:=\C^{l+m}\setminus\calA$ of the
arrangement. 
A master function complete intersection is a pure subscheme of $M_\calA$ which may be
defined by a system
 \[
   p(y)^{\beta_1}\ =\ p(y)^{\beta_2}\ =\ \dotsb\ =\ p(y)^{\beta_l}\ =\ 1
 \]
of master functions.

We describe a correspondence between systems of polynomials defining complete intersections
and systems of master functions defining  complete intersections that we call Gale
duality, as the exponent vectors of the monomials in the polynomials and 
the weights of the master functions annihilate each other.
There is also a second linear algebraic duality between the degree 1 polynomials $p_i$ and
linear forms defining the polynomials $f_i$.
Our main result is that the schemes defined by a pair of Gale dual systems are 
isomorphic.
This follows from the simple geometric observation that a complete intersection with
support $\calW$ is a linear section of the torus in an appropriate projective
embedding, and that in turn is a torus section of a linear embedding of a hyperplane
complement. 
We explain this geometry in Section~1.

In Section~2 we describe this duality concretely in terms of systems of polynomials and
systems of master functions, for this concrete version
is how it has been used. 

The value of this duality is that it allows us to transfer results about solutions to
polynomial systems to results about solutions to master functions and vice versa.
The version of this valid for positive real-number solutions was used to
give a new upper bound on the number of positive solutions of a zero-dimensional
complete 
intersection of fewnomials~\cite{BS07}, to give a continuation algorithm for
finding all 
real solutions to such a system without also computing all complex
solutions~\cite{BaSo}, 
and to give a new upper bound on the sum of the Betti numbers of a fewnomial
hypersurface~\cite{BRS07b}. 
The version valid for the real numbers leads to a surprising upper bound for the number of
real solutions to a system of fewnomials with primitive exponents~\cite{BaBiSo}.
In Section~3, we offer two results about master function complete intersections that
follow from well-known results about polynomial systems.
The first is an analog of  Kouchnirenko's bound~\cite{BKK} for the number of points in a
zero-dimensional 
master function complete intersection and the other is a formula for the Euler
characteristic of a master function complete intersection.

Another application is afforded by tropical geometry~\cite{R-GST}.
Each subvariety in the torus $(\C^\times)^{m+n}$ has an associated tropical 
variety, which is a fan in $\R^{m+n}$.
Gale duality allows us to associate certain tropical varieties to master function complete
intersections in the complement of a hyperplane arrangement.
We believe it is an interesting problem to extend this to arbitrary subvarieties of the
hyperplane complement defined by master functions.

%
%
%
\section{The geometry of Gale duality}
Let $l$, $m$, and $n$ be nonnegative integers with $l,n>0$.
We recall the standard geometric 
formulation of a system of Laurent polynomial in terms of toric varieties, 
then the less familiar geometry of systems of master functions, and then
deduce the geometric version of Gale duality.

%
%
\subsection{Sparse polynomial systems}

An integer vector $w=(a_1,\dotsc,a_{m+n})\in\Z^{m+n}$ is the exponent vector of a monomial
\[
    x^w\ :=\ x_1^{a_1}x_2^{a_2}\dotsb x_{m+n}^{a_{m+n}}\,,
\]
which is a function on the torus $(\C^\times)^{m+n}$.
Let $\calW=\{w_0,w_1,\dotsc,w_{l+m+n}\}\subset\Z^{m+n}$ be a set of exponent vectors. 
A (Laurent) polynomial $f$ with support $\calW$ is a linear
combination of monomials with exponents in $\calW$,
 \begin{equation}\label{Eq:poly}
   f(x)\ :=\ \sum_{i=0}^{l+m+n} c_i x^{w_i}\qquad\mbox{\rm where}\quad c_i\in\C\,.
 \end{equation}

A \Blue{{\it complete intersection with support $\calW$}} is a subscheme of
$(\C^\times)^{m{+}n}$ of pure dimension $n$ which may be defined by a system
 \begin{equation}\label{Eq:polynomialSystem}
   f_1(x_1,\dotsc,x_{m+n})\ =\ \dotsb\ =\  f_n(x_1,\dotsc,x_{m+n})\ =\ 0
 \end{equation}
of polynomials with support $\calW$.
Since multiplying a polynomial $f$ by a monomial does not change its zero scheme
in $(\C^\times)^{m+n}$,
we will always assume that $w_0=0$ so that our polynomials have a constant term. 

Consider the homomorphism of algebraic groups
 \begin{eqnarray*}
   \varphi_{\calW}\ \colon\ (\C^\times)^{m+n} &\longrightarrow&
     \{1\}\times(\C^\times)^{l+m+n}\ \subset\ \P^{l+m+n}\\
     x&\longmapsto& (1,x^{w_1},\dotsc,x^{w_{l+m+n}})\,.
 \end{eqnarray*}
This map $\varphi_\calW$ is dual to the homomorphism of free abelian 
groups $\Z^{l+m+n}\xrightarrow{\;\iota_\calW\;}\Z^{m{+}n}$ which maps the $i$th basis
element of $\Z^{l+m+n}$ to $w_i$. 
Write $\Z\calW$ for the image, which is the free abelian subgroup generated by $\calW$. 

The kernel of $\varphi_\calW$ is the dual $\mbox{Hom}(C\calW,\C^\times)$  
of the cokernel $C\calW:=\Z^{m+n}/\Z\calW$ of the map $\iota_\calW$.
The vector configuration $\calW$ is \Blue{{\it primitive}} when 
$\Z\calW=\Z^{m+n}$, which is equivalent to the map $\varphi_\calW$ being a closed immersion.

If we let $[z_0,z_1,\dotsc,z_{l+m+n}]$ be coordinates for $\P^{l+m+n}$, then the
polynomial $f$~\eqref{Eq:poly} equals $\varphi^*_\calW(\Lambda)$, where 
$\Lambda$ is the linear form on $\P^{l+m+n}$,
 \[
   \Lambda(z)\ =\ \sum_{i=0}^{l+m+n} c_i z_i\ .
 \]

In this way, polynomials on $(\C^\times)^{m+n}$ with support $\calW$ are pullbacks of linear
forms on $\P^{l+m+n}$.
A system~\eqref{Eq:polynomialSystem} of such polynomials defines the subscheme
$\varphi^*_\calW(L)$, where $L\subset\P^{l+m+n}$ is the linear space cut out by the forms
corresponding to the polynomials $f_i$.
An intersection $L\cap \varphi_\calW((\C^\times)^{m+n})$ is \Blue{{\it proper}} if 
its codimension equals the sum of the codimensions of $L$ and of
$\varphi_\calW((\C^\times)^{m+n})$ in $\P^{l+m+n}$. 
The following well-known proposition describes this correspondence.

\begin{proposition}\label{P:PCI}
  Every complete intersection with support $\calW$ is the pullback
  along $\varphi_\calW$ of a proper intersection of $\varphi_\calW((\C^\times)^{m+n})$ with
  a linear space $L$, and any such pullback is a complete intersection with support $\calW$.

  When $\calW$ is primitive, the map $\varphi_\calW$ is a scheme-theoretic isomorphism 
  between a complete intersection with support $\calW$ and the corresponding 
  proper intersection.
\end{proposition}

%
%
\subsection{Master functions}

Let $p_1(y),\dotsc,p_{l+m+n}(y)$ be pairwise nonproportional degree 1 polynomials on
$\C^{l+m}$. 
Their product $\prod_i p_i(y)=0$ defines a hyperplane arrangement $\calA$.
Let $\beta\in\Z^{l+m+n}$ be an integer vector, called a \Blue{{\it weight}} for the
arrangement $\calA$.
The corresponding monomial $p(y)^\beta$ in these polynomials is a 
\Blue{{\it master function}} for the arrangement $\calA$.
As the components of $\beta$ can be negative, its natural domain of definition is  
the complement $M_\calA$ of the hyperplane arrangement.

\begin{figure}[htb]
 \[   
   \begin{picture}(130,158)(0,-32)
    \put(0,0){\includegraphics[height=4.5cm]{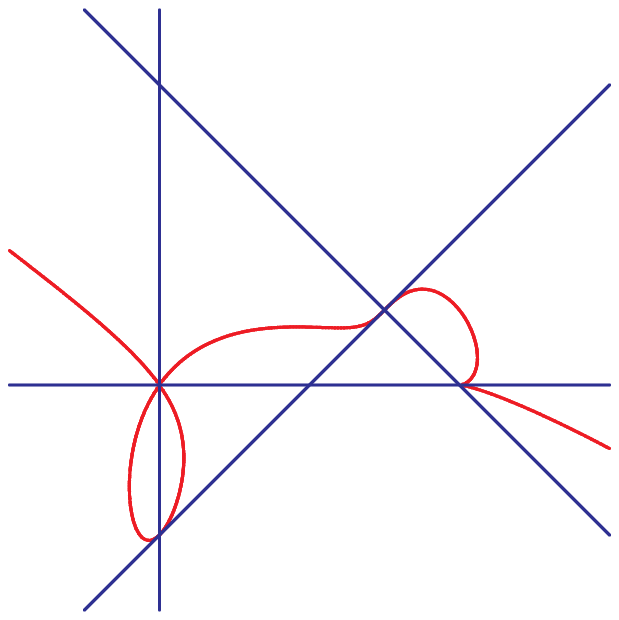}}
    \put(15,-24){${\displaystyle\frac{s^2(s+t-1)^3}{t^2(s-t-\tfrac{1}{2})}=1}$}
   \end{picture}
    \qquad\qquad
   \begin{picture}(130,158)(0,-32)
    \put(0,0){\includegraphics[height=4.5cm]{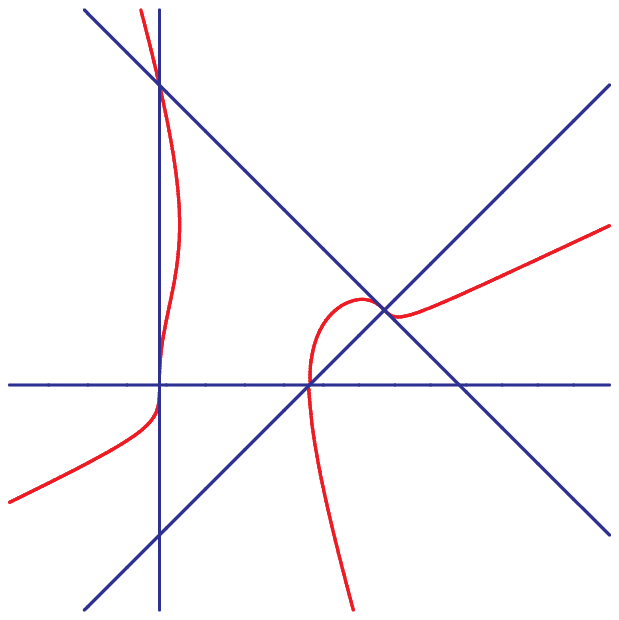}}
    \put(15,-24){${\displaystyle\frac{s(s-t-\tfrac{1}{2})^3}{t^3(s+t-1)}=1}$}
   \end{picture}
 \]
\caption{Master function curves.}\label{Fig:MF}
\end{figure}

Figure~\ref{Fig:MF} shows two curves defined by master functions in 
the complement of the arrangement
$st(s-t-\tfrac{1}{2})(s+t-1)=0$.

A \Blue{{\it master function complete intersection}} in $M_\calA$ with weights
$\calB=\{\beta_1,\dotsc,\beta_l\}$ is a subscheme of $M_\calA$ of pure dimension $m$ which
may be defined by a system 
 \begin{equation}\label{Eq:MasterSystem} 
   p(y)^{\beta_1}\ =\ p(y)^{\beta_2}\ =\ \dotsb\ =\ 
   p(y)^{\beta_l}\ =\ 1\ 
 \end{equation}
of master functions.
The weights $\calB\subset\Z^{l+m+n}$ are necessarily linearly independent.
The weights are \Blue{{\it primitive}} 
if $\Z\calB=\Q\calB\cap \Z^{l+m+n}$, so that they generated a saturated subgroup.
Linear independence and primitivity are equivalent to the subgroup
$z^{\beta_1}=\dotsb=z^{\beta_l}=1$ of the torus $(\C^\times)^{l+m+n}$
having dimension $m{+}n$ and being connected. 

The polynomials $1,p_1(y),\dotsc,p_{l+m+n}(y)$ define an affine map
 \begin{eqnarray*}
   \psi_p\ \colon\ \C^{l+m} &\longrightarrow&\P^{l+m+n}\\
                           y&\longmapsto& [1,p_1(y),\dotsc,p_{l+m+n}(y)]\,.
 \end{eqnarray*}
This map is injective if and only if the arrangement $\calA$ is \Blue{{\it essential}},
which means that the space of all degree 1 polynomials in $y$ is spanned by
$\{1,p_1(y),\dotsc,p_{l+m+n}(y)\}$.

The hyperplane arrangement $\calA$ is the pullback along $\psi_p$ of the 
coordinate hyperplanes in $\P^{l+m+n}$ and its complement $M_\calA$ is the inverse image
of the torus $(\C^{\times})^{l+m+n}\subset\P^{l+m+n}$.
Here, $\P^{l+m+n}$ has coordinates $[z_0,z_1,\dotsc,z_{l+m+n}]$ with $\psi_p^*(z_i)=p_i$.

Thus the equation $p(y)^\beta=1$ is the pullback along  
$\psi_p$ of the equation $z^\beta=1$, which defines a subgroup of
$(\C^{\times})^{l+m+n}$. 
In particular, the master function complete intersection~\eqref{Eq:MasterSystem} is the
pullback along $\psi_p$ of the subgroup $\T$ of $(\C^{\times})^{l+m+n}$
defined by $z^{\beta_1}=\dotsb=z^{\beta_l}=1$.

We summarize some properties of this corrspondence between master function complete
intersections and proper intersections of a linear space and a torus. 

\begin{proposition}\label{P:MFCI}

  Every master function complete intersection in $M_\calA$ is the pullback
  along $\psi_p$ of a proper intersection of $\psi_p(\C^{l+m})$ with a
  subgroup $\T$ of\/ $(\C^{\times})^{l+m+n}$, and any such pullback is a
  master function complete intersection in $M_\calA$.

  When $\psi_p$ is injective, the map $\psi_p$ is a scheme-theoretic isomorphism 
  between a master function complete intersection in $M_\calA$ 
  and the corresponding proper intersection.
\end{proposition}

Since no polynomial $p_i$ vanishes in $M_\calA$, we may clear denominators and rewrite 
the equation $p^\beta=1$ as an equation of polynomials, or as a binomial of the form
$p^{\beta_+}-p^{\beta_-}=0$, where 
$\beta_{\pm}$ is the vector of positive entries in $\pm\beta$.
For example, the two equations in Figure~\ref{Fig:MF} becomes
the system
\[
   s^2(s+t-1)^3- t^2(s-t-\tfrac{1}{2}) \ =\ 
   s(s-t-\tfrac{1}{2})^3 - t^3(s+t-1) \ =\ 0\,.
\]

\begin{remark}
  In the system~\eqref{Eq:MasterSystem} of master functions, each master function is set
  equal to 1.
  This is no essential restriction for if
  we instead set each master function to a different non-zero constant, then
  we may scale the degree 1 polynomials $p_i$ appropriately to absorb these constants.
  This is possible as the number of polynomials $p_i$ exceeds the number of equations.
\end{remark}

%
%
%
\subsection{Gale duality}

 Propositions~\ref{P:PCI} and~\ref{P:MFCI} form the basis of our notion of Gale duality.
 Suppose that $\T\subset(\C^{\times})^{l+m+n}\subset\P^{l+m+n}$ is a connected subgroup
 of dimension $m{+}n$ and that $L\subset\P^{l+m+n}$ is a linear space of dimension
 $l{+}m$ such that  $\T\cap L$ is proper.
 Let $L_0\simeq\C^{l+m}$ be those points of $L$ with nonzero initial (zeroth) coordinate.

\begin{definition}\label{Def:GD}
Suppose that we are given 
\begin{enumerate}
 \item An isomorphism $\varphi_\calW\colon (\C^\times)^{m+n}\to \T$ and  
         equations $z^{\beta_1}=\dotsb=z^{\beta_l}=1$  defining $\T$
         as a subgroup of $(\C^\times)^{l+m+n}$.
         Necessarily, $\calW$ and $\calB=\{\beta_1,\dotsc,\beta_l\}$ are primitive.
 \item A linear isomorphism $\psi_p\colon \C^{l+m}\to L_0$ and linear forms 
 $\Lambda_1,\dotsc,\Lambda_n$ on $\P^{l+m+n}$ defining $L$.
\end{enumerate}
Let $\calA\subset\C^{l+m}$ be the pullback of the coordinate hyperplanes of \/$\P^{l+m+n}$.
We say that the polynomial system 
 \begin{equation}\label{Eq:PS}
    \varphi_\calW^*(\Lambda_1)\ =\ \dotsb\ =\ \varphi_\calW^*(\Lambda_n)\ =\ 0
 \end{equation}
in $(\C^\times)^{m+n}$ is \Blue{{\it Gale dual}} to the system of master functions
 \begin{equation}\label{Eq:MFS}
   \psi_p^*(z^{\beta_1})\ =\ \dotsb\ =\ \psi_p^*(z^{\beta_l})\ =\ 1
 \end{equation}
in $M_\calA$ and vice-versa.
\end{definition}

This definition contains two different linear algebra dualities.
The weights $\calB$ form a $\Z$-basis for the free abelian group
of integer linear relations among the nonzero exponent vectors 
of $\calW$.
Similarly, the linear forms $\{\Lambda_1,\dotsc,\Lambda_n\}$ form a basis for the 
space of linear relations among the coordinate functions $\{1,p_1,\dotsc,p_{l+m+n}\}$
defining the map $\psi_p$.

The following is immediate.

\begin{theorem}\label{Th:GD}
  A pair of Gale dual systems~$\eqref{Eq:PS}$ and~$\eqref{Eq:MFS}$ define isomorphic
  schemes. 
\end{theorem}

%
%
%
\section{The algebra of Gale duality}

We give an explicit algorithmic version of Gale duality.
Let $\calW=\{0,w_1,\dotsc,w_{l+m+n}\}\subset\Z^{m+n}$ be a primitive collection of integer
vectors and suppose that 
 \begin{equation}\label{Eq:PCI}
   f_1(x_1,\dotsc,x_{m+n})\ =\  
    \dotsb\ =\  f_n(x_1,\dotsc,x_{m+n})\ =\ 0
 \end{equation}
defines a complete intersection with support $\calW$ in the torus $(\C^\times)^{m+n}$.
Then the polynomials $f_i$ are linearly independent. 
We may reorder the exponent vectors so that the coefficients
of $x^{w_1},\dotsc,x^{w_n}$ in~\eqref{Eq:PCI} form an invertible matrix
and then transform~\eqref{Eq:PCI} into an equivalent system
where the coefficients of $x^{w_1},\dotsc,x^{w_n}$ form a diagonal matrix.
 \begin{equation}\label{Eq:Diag_system}
  \begin{array}{rcl}
  x^{w_1} &=& g_1(x)\ =:\ p_1(x^{w_{n+1}},\dotsc,x^{w_{l+m+n}})\\\rule{0pt}{14pt}
          &\vdots&\\\rule{0pt}{14pt}
  x^{w_n} &=& g_n(x)\ =:\ p_n(x^{w_{n+1}},\dotsc,x^{w_{l+m+n}})
  \end{array}
 \end{equation}

Here, for each $i=1,\dotsc,n$, $g_i(x)$ 
is a polynomial with support $\{0,\, w_{n+1},\dotsc,w_{l+m+n}\}$ which is a 
degree 1 polynomial function $p_i(x^{w_{n+1}},\dotsc,x^{w_{l+m+n}})$ in the given $l{+}m$ arguments.
For $i=n{+}1,\dotsc,l{+}m{+}n$, set 
$p_i(x^{w_{n+1}},\dotsc,x^{w_{l+m+n}}):= x^{w_i}$.

An integer linear relation among the exponent vectors in $\calW$,
\[
  b_1 w_1\;+\; b_2w_2\;+\;\dotsb\;+\;
   b_{l+m+n}w_{l+m+n}\ =\ 0\,,
\]
is equivalent to the monomial identity
\[
   (x^{w_1})^{b_1}\cdot (x^{w_2})^{b_2}\dotsb(x^{w_{l+m+n}})^{b_{l+m+n}}
    \ =\ 1\,,
\]
which gives the consequence of the system~\eqref{Eq:Diag_system}
\[
   \bigl(p_1(x^{w_{n+1}},\dotsc,x^{w_{l+m+n}})\bigr)^{b_1}\ \dotsb\ 
   \bigl(p_{l+m+n}(x^{w_{n+1}},\dotsc,x^{w_{l+m+n}})\bigr)^{b_{l+m+n}}
    \ =\ 1\,.
\]

Define $y_1,\dotsc,y_{l+m}$ to be new variables which are coordinates for 
$\C^{l+m}$.
The degree 1 polynomials $p_i(y_1,\dotsc,y_{l+m})$ define a hyperplane arrangement
$\calA$ in $\C^{l+m}$.
Note that $\calA$ is essential since it contains all the coordinate hyperplanes of $\C^{l+m}$.
Let  $\calB:=\{\beta_1,\dotsc,\beta_m\}\subset\Z^{l+m+n}$ be a basis for the 
$\Z$-module of integer linear relations among the nonzero vectors in $\calW$.
These weights $\calB$ define a system of master functions 
 \begin{equation}\label{Eq:MFCI}
   p(y)^{\beta_1}\ =\ p(y)^{\beta_2}\ =\ \dotsb\ =\ 
  p(y)^{\beta_l}\ =\ 1
 \end{equation}
in the complement $M_\calA:=\C^{l+m}\setminus \calA$ of the arrangement.

\begin{theorem}\label{Thm:Alg-GD}
  The system of polynomials~\eqref{Eq:PCI} in $(\C^\times)^{m+n}$
  and the system of master functions~\eqref{Eq:MFCI} in  
  $M_\calA$ define isomorphic complete intersections.
\end{theorem}

\begin{proof}
 Condition (1) in Definition~\ref{Def:GD} holds as $\calW$ and $\calB$ are both primitive 
 and annihilate each other.
 The linear forms $\Lambda_i$ that pull back along $\varphi_\calW$ to 
 define the system~\eqref{Eq:Diag_system} are
\[
   \Lambda_i(z)\ =\ z_i\ -\ p_i(z_{n{+}1},\dotsc,z_{l+m+n})\,,
\]
 which shows that condition (2) holds, and so the statement follows from
 Theorem~\ref{Th:GD}.
\end{proof}

\begin{example}\label{Ex:main}
  Suppose that we have the system of polynomial equations
 \begin{equation}\label{Eq:example}
   \begin{array}{rcl}
     2x^4y^{-1} - 3x^3y^2 - 4x^4y + xy^2 -\tfrac{1}{2} &=&0\\\rule{0pt}{14pt}
      x^3y^2 + 2x^4y - xy^2 -\tfrac{1}{2}&=&0
   \end{array}
 \end{equation}
 in the torus $(\C^\times)^2$.
 Here $n=l=2$ and $m=0$.
 We may diagonalize this to obtain 
 \begin{eqnarray*} 
   x^3y^2 &=& x^4y^{-1}-x^4y-\tfrac{1}{2}\,,\\
   xy^2   &=& x^4y^{-1}+x^4y -1\,.
 \end{eqnarray*}
 Thus the system has the form
 $\varphi^{*}_\calW(\Lambda_1)=\varphi^{*}_\calW(\Lambda_2)=0$,
 where
\[
   \begin{array}{rcl}
   \Lambda_1(z) &=& z_1-(z_3-z_4-\tfrac{1}{2})\,,\\\rule{0pt}{14pt}
   \Lambda_2(z) &=& z_2-(z_3+z_4-1)\,,    \qquad\mbox{and}\qquad\\\rule{0pt}{14pt}
    \varphi_\calW&\colon&(x,y)\ \longmapsto\  
       (x^3y^2,xy^2,x^4y^{-1},x^4y)\ =\ (z_1,z_2,z_3,z_4)\,.
   \end{array}
\]
 These exponents $\calW$ are primitive.

 Let $s,t$ be new variables and set
 \begin{eqnarray*}
   p_1\ :=\ s-t-\tfrac{1}{2} &\qquad& p_3\ :=\ s\\
   p_2\ :=\ s+t-1&& p_4\ :=\ t
 \end{eqnarray*}
 Then $\psi_p\colon (s,t)\mapsto (p_1,p_2,p_3,p_4)$ parametrizes the
 common zeroes of $\Lambda_1$ and $\Lambda_2$.

 Note that 
\[
   (x^3y^2)^{-1}(xy^2)^3   (x^4y^{-1})^2(x^4y)^{-2}\ =\ 
   (x^3y^2)^3   (xy^2)^{-1}(x^4y^{-1})  (x^4y)^{-3}\ =\  1\,,
\]
 and so the weights  $(-1,3,2,-2)$ and $(3,-1,1,-3)$ annihilate $\calW$.
 These weights are primitive.
 By Theorem~\ref{Thm:Alg-GD}, the polynomial system~\eqref{Eq:example}
 in $(\C^\times)^2$ is equivalent to the system of master functions
 \begin{equation}\label{Eq:MFex}
  \begin{array}{rcl}
    s^2(s+t-1)^3\ -\ t^2(s-t-\tfrac{1}{2})&=&0\\\rule{0pt}{14pt}
    s(s-t-\tfrac{1}{2})^3\ -\ t^3(s+t-1) &=& 0\,,
  \end{array}
 \end{equation}
 in the complement of the hyperplane arrangement
 $st(s+t-1)(s-t-\tfrac{1}{2})=0$.
 We display these two systems in Figure~\ref{Fig:one}, drawing also the
 excluded hyperplanes (lines).
\begin{figure}[htb]
\[
   \begin{picture}(160,150)
    \put(0,0){\includegraphics[height=5.2cm]{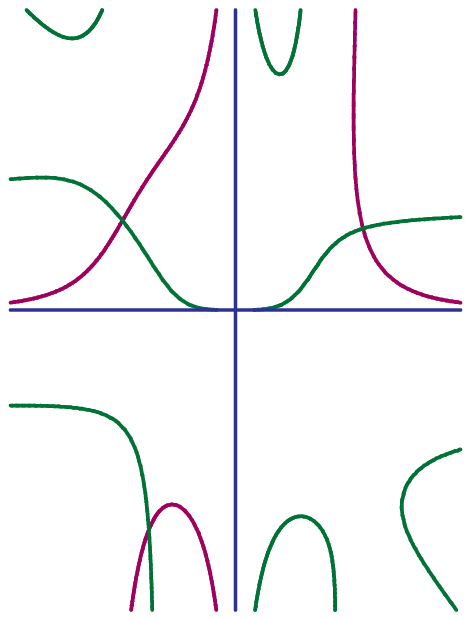}}
    \put(43,-10){$x=0$}    \put(118, 74){$y=0$}
   \end{picture}
    \qquad\qquad
   \begin{picture}(195,150)
    \put(0,0){\includegraphics[height=5.2cm]{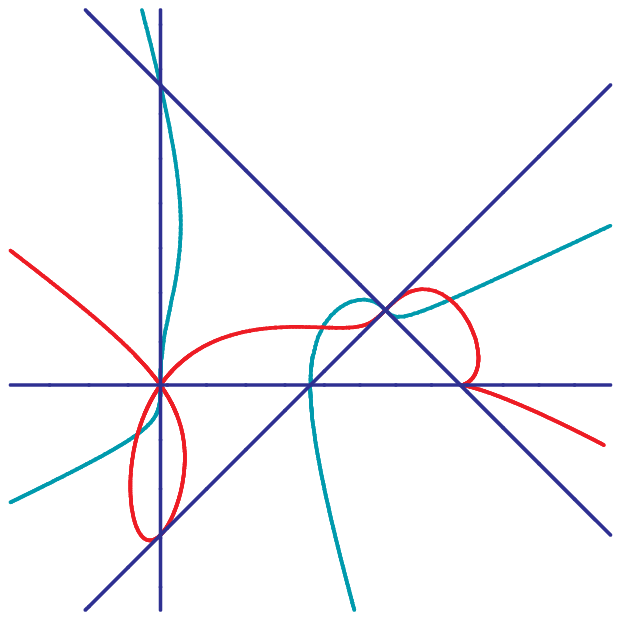}}
    \put(25,-10){$s=0$}   \put(154,54){$t=0$}
    \put(120,5){$s+t-1=0$}  \put(110,140){$s-t-\frac{1}{2}=0$}
   \end{picture}
\]
\caption{The polynomial system~\eqref{Eq:example} and the 
      system of master functions~\eqref{Eq:MFex}.}\label{Fig:one}
\end{figure}

\end{example}

We remark that although we have two curves in the polynomial system and two curves in the
system of master functions, the individual curves are unrelated.
Theorem~\ref{Thm:Alg-GD} merely asserts an isomorphism between the zero-dimensional schemes 
in the torus $(\C^\times)^2$ and in the hyperplane complement $M_\calA$ defined by
each pair of curves.
%
%
%
\section{Some consequences of Gale duality}

Theorems~\ref{Th:GD} and~\ref{Thm:Alg-GD}, which assert isomorphisms of schemes,
hold if $\C$ is replaced by any other field, and even other algebraic objects.
In particular, Gale duality holds for the real numbers.
Example~\ref{Ex:main} illustrates this fact.
The real zero-dimensional schemes defined by~\eqref{Eq:example} in $(\R^\times)^2$ and 
by~\eqref{Eq:MFex} in $M_\calA(\R)$ each consist of 3 reduced points with residue field
$\R$ (which we see in Figure~\ref{Fig:one}) and 7 reduced points with residue field $\C$.
If we  only consider real-number solutions, that is, analytic subschemes of
$(\R^\times)^{m+n}$ and of $M_\calA(\R)$, then we may relax the requirement in
Theorem~\ref{Th:GD} that $\calW$ and $\calB$ are primitive to the condition that 
they generate subgroups of odd index in their saturations.

Gale duality also holds for $\R_+$, the positive real numbers and for
$M^+_\calA$, the positive chamber of the complement $M_\calA(\R)$ of oriented hyperplanes
$\calA$. 
In this guise it is Theorem~2.2 of~\cite{BS07}.
There, positivity allows $\calW$ and $\calB$ to have real-number
components. \smallskip

Gale duality allows us to use knowledge about polynomial systems to deduce results about
systems of master functions, and vice-versa.
In fact, this is how it arose.
It was used implicitly~\cite{BBS,B07} and explicitly~\cite{BS07} to give new upper
bounds on the number of positive solutions to a system of fewnomial equations. 
In~\cite{BS07}, the bound
 \begin{equation}\label{Eq:Fewnomial}
   \frac{e^2+3}{4} 2^{\binom{l}{2}}n^l
 \end{equation}
was given for the number of solutions to a 0-dimensional
master function complete intersection ($m=0$) in the positive 
chamber $M^+_\calA$, where $\calA$ consists of $l{+}n$ oriented hyperplanes in $\R^l$.
By Gale duality for $\R_+$, we obtain the new fewnomial bound of~\cite{BS07}:
A system of $n$ polynomials in $n$ variables having a total of $n{+}l{+}1$ distinct monomials 
has at most~\eqref{Eq:Fewnomial} nondegenerate solutions in the positive orthant $\R_+^n$.

The proof in~\cite{BS07} leads to a path continuation algorithm~\cite{BaSo} to find 
nondegenerate solutions in $M^+_\calA$ to systems of master functions.
Its novelty is that, unlike traditional continuation algorithms for 
solving systems of algebraic equations~\cite{SW05}, it only follows
real solutions. 
Its complexity depends upon the dimension $l$ and the fewnomial bound~\eqref{Eq:Fewnomial},
and not on the number of complex solutions to the system of master functions.
That algorithm easily extends to find all nondegenerate solutions in the hyperplane complement
$M_\calA(\R)$, and through Gale duality it gives a new continuation algorithm for 
all nondegenerate real solutions to a system of polynomial equations.
Moreover, the ideas underlying the algorithm lead to a generalization of~\cite{BS07},
giving the bound
\[
   \frac{e^4+3}{4} 2^{\binom{l}{2}}n^l
\]
for the number of nondegenerate real solutions to a master function complete 
intersection and thus a bound for the number of nonzero nondegenerate 
real solutions to a system of $n$ polynomials in
$n$ variables having a total of $n{+}l{+}1$ distinct monomials~\cite{BaBiSo}.

These new fewnomial bounds are used to bound the number of connected
components~\cite{BRS07a} and the sum of the Betti numbers~\cite{BRS07b}
of a fewnomial hypersurface.
For example, Theorem~1 of~\cite{BRS07b} states that
the sum of the Betti numbers of a hypersurface in $\R^{m+1}_+$
defined by a polynomial with $l{+}m{+}1{+}1$ monomial terms is bounded by 
 \begin{equation}\label{Eq:betti}
  \frac{e^2+3}{4} 2^{\binom{l}{2}} (m+1)^l\cdot 2^{m+1}\,.
 \end{equation}
By Gale duality, this gives a following bound for certain complete intersections of 
master functions.

\begin{corollary}
  Let $\calA$ be an arrangement in $\R^{l+m}$ consisting of\/ 
  $l{+}m{+}1$ hyperplanes.
  Let $M^+_\calA$ be a chamber of the hyperplane complement $M_\calA$.
  Then the sum of the Betti numbers of a smooth codimension $l$ complete 
  intersection defined by master functions in $M^+_\calA$ is at most~$\eqref{Eq:betti}$.
\end{corollary}
\begin{proof}
 The arrangement $\calA$ has one hyperplane more than the dimension of the ambient space
 $\R^{l+m}$.
 Therefore $n=1$, and the Gale dual complete intersection is a smooth hypersurface
 in $\R^{m+1}_+$. 
\end{proof}

Another interesting class of applications of Gale duality is to transfer
results about polynomial systems (which have been extensively studied)
to systems of master functions, which have not yet attracted 
much attention.

Let $\calB:=\{\beta_1,\dotsc,\beta_l\}$ be linearly independent elements of 
$\Z^{l+m+n}$ which are primitive (that is, $\Z\calB=\Q\calB\cap \Z^{l+m+n}$).
Then the quotient
\[
   \Z^{l+m+n}/\Z\calB
\]
is a free abelian group of rank $m{+}n$ that we identify with 
$\Z^{m+n}$.
For each $i=1,\dotsc,l{+}m{+}n$, let $w_i\in\Z^{m+n}$ be the image of the 
$i$th standard unit vector in $\Z^{l+m+n}$.
These generate $\Z^{m+n}$ and so $\Blue{\calW}:=\{0,w_1,\dotsc,w_{l+m+n}\}$, is 
primitive.
Let
\[
  \Blue{\Delta_\calB}\ :=\ \mbox{conv}(0,w_1,\dotsc,w_{l+m+n})
\]
be the convex hull of $\calW$.

Our first application is a Kouchnirenko Theorem~\cite{BKK} for zero-dimensional
($m=0$) complete intersections of master functions.

\begin{corollary}[Kouchnirenko's Theorem for master functions]\label{Co:KTh}
  Let $p_1,\dotsc,p_{l+n}$ be degree $1$ polynomials which define an essential
  arrangement $\calA$ of $l+n$ hyperplanes in $\C^l$.
  Then the system of master functions
 \begin{equation}\label{Eq:THM_MF}
   p^{\beta_1}\ =\ p^{\beta_2}\ =\ \dotsb\ =\ p^{\beta_l}\ =\ 1\,,
 \end{equation}
 in the hyperplane complement $M_\calA$ has at most
 \begin{equation}\label{Eq:Kouch_num}
   n!\; \mbox{\rm volume}(\Delta_\calB)
 \end{equation}
 isolated solutions, counted with multiplicity.
 When the polynomials $p_1,\dotsc,p_{l+n}$ are general, the
 system~$\eqref{Eq:THM_MF}$ has exactly~$\eqref{Eq:Kouch_num}$ solutions.
\end{corollary}

\begin{example}\label{Ex:Kouchnirenko}
 For example, the master functions of Figure~\ref{Fig:MF} and of the
 system~\eqref{Eq:MFex} have weights
\[
  \calB\ :=\ \{ (-1,3,2,-2),\ (3,-1,1,-3)\}\,.
\]
 These are primitive, so $\Z^4/\Z\calB\simeq \Z^2$.
 This isomorphism is realized by sending the standard basis vectors of $\Z^4$ to the
 columns of the matrix
\[
   \calW\ =\ \left[\begin{array}{rrrr}3&1&4&4\\2&2&-1&1\end{array}\right]\ .
\]
 The convex hull of these columns and the origin is the integer pentagon
 in $\R^2$
\[
   \includegraphics{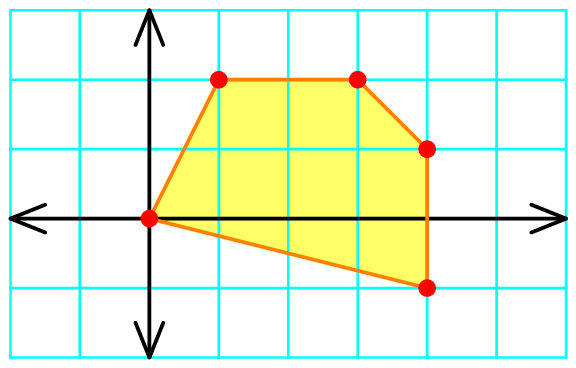}
\]
 which has area $17/2$.
 Thus by Corollary~\ref{Co:KTh}, a system of master functions for a general arrangement of
 4 lines in $\C^2$ with weights $\calB$ will have $17$ solutions in the complement of 
 of the arrangement.
 Indeed, the system~\eqref{Eq:MFex} of master functions has 17 solutions in the complement
 of the line arrangement shown in Figure~\ref{Fig:one}.
 Similarly, the system 
\[
   \frac{(2x-3y)^2(4x+y-7)^3}{(1+x-3y)^2(x-7y-2)}\ =\ 
    \frac{(2x-3y)(x-7y-2)^3}{(1+x-3y)^3(4x+y-7)}\ =\ 1\,
\]
 also has 17 solutions in the complement of its line arrangement.
 These claims of 17 solutions are readily checked by computer.
\end{example}

\begin{proof}[Proof of Corollary~$\ref{Co:KTh}$]
 The polynomials $1,p_1,\dotsc,p_{l+n}$ parametrize a codimension $n$
 plane $\Lambda\subset\P^{l+n}$, and the system~\eqref{Eq:THM_MF} 
 defines the intersection
 \begin{equation}\label{Eq:int}
   \varphi_\calW( (\C^\times)^n) \cap \Lambda\,.
 \end{equation}
 Since $\Lambda$ has codimension $n$, this is a complete intersection
 with support $\calW$.
 The first statement follows by Kouchnirenko's Theorem~\cite{BKK}.

 For the second statement, observe that if $\Lambda$ is a general 
 codimension $n$ plane, then the intersection~\eqref{Eq:int} is transverse
 and Kouchnirenko's theorem implies that it consists of exactly~\eqref{Eq:Kouch_num}
 points.
 But a general codimension $n$ plane $\Lambda$ in $\P^{l+n}$ is parametrized
 by general polynomials $1,p_1,\dotsc,p_{l+m+n}$.
\end{proof}

Khovanskii~\cite{Kh78} gave formulas for many invariants of complete intersections
in the torus, including genus, arithmetic genus, and Euler characteristic.
By Gale duality, these are formulas for invariants of master function 
complete intersections. 
Khovanskii's formulas for genera are rather involved, and we leave their formulation for
master function complete intersections as an exercise for the interested reader.
His formula for the Euler characteristic is however quite simple.
Let $\calB$ and $\Delta_\calB$ be as described before Theorem~\ref{Co:KTh}.

\begin{corollary}
  Let $p_1,\dotsc,p_{l+m+n}$ be general degree $1$ polynomials  which define an essential
  arrangement $\calA$ of $l{+}m{+}n$ hyperplanes in $\C^{l+m}$.
  The Euler characteristic of the solution set of the system of master functions
 \begin{equation}\label{Eq:THM_nMF}
   p^{\beta_1}\ =\ p^{\beta_2}\ =\ \dotsb\ =\ p^{\beta_l}\ =\ 1\,,
 \end{equation}
 is 
 \begin{equation}\label{Eq:Kouch_nnum}
   (-1)^m \binom{m{+}n{-}1}{n-1} \cdot (m{+}n)!\;\mbox{\rm volume}(\Delta_\calB)\,.
 \end{equation}
\end{corollary}

\noindent{\it Proof.}
  We compute the Euler characteristic of the complete intersection in a torus defined by a
  system of polynomials Gale dual to the master functions in~\eqref{Eq:THM_nMF}.
  Khovanskii~\cite[Section 3, Theorem 1]{Kh78} shows that the Euler 
  characteristic of a transverse intersection 
  of hypersurfaces $X_1,\dotsc,X_n$ in $(\C^\times)^{m+n}$ is 
 \begin{equation}\label{Eq:Khovanskii}
   \Bigl( \prod_{i=1}^n \frac{D_i}{1+D_i} \Bigr)\ \cap\ [(\C^\times)^{m+n}]\,,
 \end{equation}
  where $D_i$ is the divisor class of $X_i$ and $[(\C^\times)^{m+n}]$ 
  is the fundamental class of the torus $(\C^\times)^{m+n}$.
  This is computed in the Chow ring~\cite{FuSt} of any toric variety where it makes
  sense.  

  In our application of his result, the divisors are equal, say to $D$,
  and we have 
\[
   D^{m+n}\cap  [(\C^\times)^{m+n}]\ =\ (m{+}n)!\;\mbox{volume}(\Delta_\calB)\,,
\]
 by Kouchnirenko's Theorem.
 Thus Khovanskii's formula~\eqref{Eq:Khovanskii} becomes
 \begin{eqnarray*}
   \Bigl( \frac{D}{1+D}\Bigr)^n \cap [(\C^\times)^{m+n}] &=&
   \Bigl( D\sum_{j\geq 0} (-1)^j D^j \Bigr)^n \cap [(\C^\times)^{m+n}]\\
    &=& \Bigl((-1)^m  \sum_{j_1+\dotsb+j_n=m} 1 \Bigr)
         \cdot D^{m+n} \cap [(\C^\times)^{m+n}]\\
    &=&    (-1)^m \binom{m+n-1}{n-1} \cdot (m{+}n)!\;\mbox{\rm volume}(\Delta_\calB)\,.
     \makebox[1pt][l]{\hspace{38pt}$\Box$}\vspace{20pt}
 \end{eqnarray*}


\def\cprime{$'$}
\providecommand{\bysame}{\leavevmode\hbox to3em{\hrulefill}\thinspace}
\providecommand{\MR}{\relax\ifhmode\unskip\space\fi MR }
\providecommand{\MRhref}[2]{%
  \href{http://www.ams.org/mathscinet-getitem?mr=#1}{#2}
}
\providecommand{\href}[2]{#2}

\end{document}